\documentclass[12pt]{article}
\usepackage{amsmath,amssymb}
\usepackage{rotating}
\usepackage{xypic}
\usepackage{amssymb}
\usepackage{amsmath,amscd}
\usepackage{pst-node,pstricks,multido,pst-plot,pst-text,pst-3d}%
\usepackage{graphicx}
\usepackage{amsfonts}
\newtheorem{example}{Example}[section]

\newtheorem{remark}[example]{Remark}
\newtheorem{theorem}[example]{Theorem}

\newtheorem{corollary}[example]{Corollary}

\newtheorem{proposition}[example]{Proposition}

\newtheorem{lemma}[example]{Lemma}

\def\S{{\mathfrak  S}}
\def\cal#1{{\mathfrak #1}}

\def\<{\langle}
\def\>{\rangle}

\def\Z{{\mathbb Z}}
\def\N{{\mathbb N}}
\def\Y{{\mathbb Y}}

\def\CT{{\rm C.T.}}

\def\goth{\mathfrak}

\def\ashuff#1#2#3{
\kern 1pt \vrule height#1 \overline{\vrule height#3 width 0pt
\hskip#2} \rule{.3pt}{#1}\overline{\vrule height#3 width 0pt
\hskip#2} \rule{.3pt}{#1} \kern 1pt }

\def\X{{\mathbb X}}\def\Y{{\mathbb Y}}

\def\sign{\mbox{sign}}

\begin{document}

\title{Macdonald polynomials at $t=q^k$}
\author{Jean-Gabriel Luque}
 \maketitle
 \begin{abstract}
We investigate the homogeneous symmetric Macdonald polynomials $P_\lambda(\X;q,t)$ for the specialization $t=q^k$.
We show  an identity relying
 the polynomials
$P_\lambda(\X;q,q^k)$ and $P_\lambda\left(\frac{1-q}{1-q^k}\X;q,q^k\right)$.
 As a consequence, we describe an operator whose eigenvalues characterize the polynomials $P_\lambda(\X;q,q^k)$.
 \end{abstract}
\section{Introduction}
Macdonald polynomials  are $(q,t)$-deformations of Schur functions which play an important r\^ole in the representation theory
of the double affine Hecke algebra \cite{Lasc2,MacdoHecke} since they are the eigenfunctions of the Cherednik elements.
The polynomials considered here are the homogeneous symmetric Macdonald polynomials $P_\lambda(\X;q,t)$ and are the eigenfunctions of the
Sekiguchi-Debiard operator.
 For $(q,t)$ generic, these polynomials are completely
characterized by their eigenvalues, since the dimensions of the eigenspaces is $1$. It is no longer the case when $t$
is specialized to
a rational power of $q$. Hence, it is more convenient to characterize the Macdonald (homogeneous symmetric)
polynomials by orthogonality  ({\it w.r.t.}  a ($q,t$)-deformation of the usual scalar product on symmetric functions) and
by some conditions on their dominant monomials (see {\it e.g.} \cite{Macdo}).
 In this paper, we consider the specialization $t=q^k$ where $k$ is a strictly positive
integer.
One of our motivations is to generalize an identity of \cite{BBL}, which shows that even powers of the discriminant
are rectangular Jack polynomials. Here, we show  that this property follows from deeper relations between
 the Macdonald polynomials
$P_\lambda(\X;q,q^k)$ and $P_\lambda\left(\frac{1-q}{1-q^k}\X;q,q^k\right)$ (in the $\lambda$-ring notation).
This result is interesting in the context of the quantum fractional Hall effect\cite{Lau},
 since it implies  properties of the expansion of the powers of the discriminant in the Schur basis
\cite{DGIL,KTW,STW}.
It implies also that the Macdonald polynomials (for $t=q^k$) are
 characterized by the eigenvalues of an operator $\cal M$ whose eigenspaces are of dimension  $1$
described in terms of isobaric divided differences.

The paper is organized as follow. After recalling notations and background (Section \ref{SNot}) for
 Macdonald polynomials,
we give, in Section \ref{Ssubsti}, some properties of the operator which substitutes a complete function to each power
of a letter. These properties allow to show our main result in Section \ref{Smain}
which is an identity relying the polynomial
$P_\lambda(\X;q,q^k)$ and $P_\lambda\left(\frac{1-q}{1-q^k}\X;q,q^k\right)$.
 As a consequence, we  describe  (Section \ref{Seigen})
an operator ${\cal M}$ whose eigenvalues characterize the Macdonald polynomials $P_\lambda(\X;q,q^k)$.
 Finally, in Section \ref{SChered}, we give an
expression of ${\cal M}$ in terms of Cherednik elements.

\section{Notations and background\label{SNot}}
Consider an alphabet $\X$ potentially infinite. We will use the
notations of \cite{Lasc} for the generating function $\sigma_z(\X)$
of the complete homogeneous functions $S^p(\X)$,
\[
\sigma_z(\X)=\sum_iS^i(\X)z^i=\prod_i\frac1{1-xz}.
\]
The algebra $Sym$ of symmetric function has a structure of
$\lambda$-ring \cite{Lasc}. We recall that the sum of two alphabets $\X+\Y$ is defined by
\[\sigma_z(\X+\Y)=\sigma_z(\X)\sigma_z(\Y)=\sum_{i}S^i(\X+\Y)z^i.\]
 In particular, if $\X=\Y$ one has $\sigma_z(2\X)=\sigma_z(\X)^2$.
 This definition is extended for any complex number $\alpha$ by
 $\sigma_z(\alpha\X)=\sigma_z(\X)^\alpha$.
 For example, the generating series of the elementary function is
 \[\begin{array}{rcl}\lambda_z(\X)&:=&\sum\Lambda_i(\X)z^i=\prod_x(1+xz)\\
 &=&\sigma_{-z}(-\X)=\sum_i(-1)^iS^i(-\X)z^i.\end{array}\]
 The complete functions of the product of two alphabets $\X\Y$ are
given by the Cauchy kernel
\[
K(\X,\Y):=\sigma_1(\X\Y)=\sum_iS^i(\X\Y)=\prod_{x\in\X}\prod_{y\in\Y}{1\over1-xyt}=\sum_\lambda S_\lambda(\X)S_\lambda(\Y),
\]
where $S_\lambda$ denotes, as in \cite{Lasc}, a Schur function. More generally, one has
\[
K(\X,\Y)=\sum_\lambda A_\lambda(\X)B_\lambda(\Y)
\]
for any pair of basis $(A_\lambda)_\lambda$ and $(B_\lambda)_\lambda$ in duality for the usual scalar product
$\langle\,,\,\rangle$.

\subsection{Macdonald polynomials}

One considers the $(q,t)$-deformation (see {\it e.g.} \cite{Macdo}) of the usual scalar product on
symmetric functions defined for a pair of power sum functions $\Psi^\lambda$ and $\Psi^\mu$
(in the notation of \cite{Lasc}) by
\begin{equation}
\langle\Psi^\lambda,\Psi^\mu\rangle_{q,t}=
\delta_{\lambda,\mu}z_\lambda\prod_{i=1}^{l(\lambda)}{1-q^{\lambda_i}\over
1-t^{\lambda_i}},
\end{equation}
where $\delta_{\lambda,\mu}=1$ if $\lambda=\mu$ and $0$ otherwise.
The familly of Macdonald polynomials $(P_\lambda(\X;q,t))_\lambda$ is the
unique basis of symmetric functions orthogonal for $\langle\ ,\
\rangle_{q,t}$ verifying
\begin{equation}
P_\lambda(\X;q,t)=m_\lambda(\X)+\sum_{\mu\leq\lambda}u_{\lambda\mu}m_\mu(\X),
\end{equation}
where $m_\lambda$ denotes, as usual, a monomial function \cite{Lasc,Macdo}.
The reproducing kernel associated to this scalar product is
\[K_{q,t}(\X,\Y):=
\sum_\lambda\langle\Psi^\lambda,\Psi^\lambda\rangle_{q,t}^{-1}\Psi_\lambda(\X)\Psi_\lambda(\Y)
=\sigma_1\left({1-t\over 1-q}\X\Y\right)\]
see {\it e.g.} \cite{Macdo} (VI. 2). In particular, one has
\begin{equation}\label{Kqt}K_{q,t}(\X,\Y)=\sum_\lambda
P_{\lambda}(\X;q,t)Q_{\lambda}(\Y;q,t),\end{equation}
where $Q_\lambda(\X;q,t)$ is the dual  basis of $P_\lambda(\Y;q,t)$
for $\langle\, ,\, \rangle_{q,t}$,
\begin{equation}
Q_{\lambda}(\X;q,t)=\langle P_\lambda,
P_\lambda\rangle_{q,t}^{-1}P_\lambda(\X;q,t).
\end{equation}
The coefficient $b_\lambda(q,t)=\langle P_\lambda,
P_\lambda\rangle_{q,t}^{-1}$ is known to be
\begin{equation}
b_{\lambda}(q,t)=\prod_{(i,j)\in\lambda}
{1-q^{\lambda_j-i+1}t^{\lambda'_i-j}\over
1-q^{\lambda_j-i}t^{\lambda'_i-j+1}}
\end{equation}
see \cite{Macdo} VI.6. Writing
\begin{equation}\label{Kqt2K}
K_{q,t}\left(\left(1-q\over1-t\right)\X,\Y\right)=K(\X,\Y),
\end{equation}
one finds that
$\left(P_\lambda\left(\left(1-q\over1-t\right)\X;q,t\right)\right)_\lambda$
is the dual basis of $\left(Q_\lambda(\X;q,t)\right)_\lambda$ for
the usual scalar product $\langle\ ,\ \rangle$.

Note that there exists an other Kernel type formula which reads
\begin{equation}\label{lambda_t}
\lambda_1(\X\Y)=\sum_\lambda
P_{\lambda'}(\X;t,q)P_\lambda(\Y;q,t)=\sum_\lambda
Q_{\lambda'}(\X;t,q)Q_\lambda(\Y;q,t).
\end{equation}
where $\lambda'$ denotes the conjugate partition of $\lambda$.
This formula can be found in \cite{Macdo} VI.5 p 329.

From Equalities (\ref{Kqt2K}) and (\ref{Kqt}) , one has
\begin{equation}\label{-X1}
\sigma_1(\X\Y)=K_{q,t}\left({1-q\over1-t}\X,\Y\right)=
\sum_\lambda Q_{\lambda}\left({1-q\over 1-t}\X;q,t\right)
P_\lambda(\Y;q,t).
\end{equation}
Applying (\ref{lambda_t}) to
\[ \sigma_1(\X\Y)=\lambda_{-1}(-\X\Y),\]
one obtains
\begin{equation}\label{-X2}
\sigma_1(\X\Y)=\sum_{\lambda}(-1)^{|\lambda|}Q_{\lambda'}(-\X;t,q)Q_{\lambda}(\Y;q,t).
\end{equation}
Identifying the coefficient of $P_\lambda(\Y;t,q)$ in (\ref{-X1}) and (\ref{-X2}), one finds the property below.
\begin{lemma}\label{P-X2Pqt}
\begin{equation}
Q_\lambda(-\X;t,q)=(-1)^{|\lambda|}P_{\lambda'}\left({1-q\over 1-t}\X;q,t\right).
\end{equation}
\end{lemma}

Unlike the usual ($q=t=1$) scalar product, there is  no expression as
a constant term for the product $\langle\,,\,\rangle_{q,t}$ when
$\X=\{x_1,\dots,x_n\}$ is finite. But the Macdonald polynomials are
orthogonal for an other scalar product defined by
\begin{equation}\label{MacdoScalarPrime}
\langle
f,g\rangle'_{q,t;n}=\frac1{n!}\CT\{f(\X)g(\X^\vee)\Delta_{q,t}(\X)\}
\end{equation}
where $\CT$ denotes constant term {\it w.r.t.} the alphabet $\X$, $\Delta_{q,t}(\X)=\displaystyle\prod_{i\neq
j}{(x_ix_j^{-1};q)_\infty \over(tx_ix_j^{-1};q)_\infty}$,
$(a;b)_\infty=\displaystyle\prod_{i\geq 0}(1-ab^i)$ and $\X^\vee=\{x_1^{-1},\dots,x_n^{-1}\}$.
The expression of $\langle P_\lambda,Q_\lambda\rangle'_{q,t;n}$ is given by (\cite{Macdo} VI.9)
\begin{equation}\label{PQprime}
\langle P_\lambda,Q_\lambda\rangle_{q,t;n}'=\frac1{n!}\CT\{\Delta_{q,t}(\X)\}
\prod_{(i,j)\in \lambda}{1-q^{i-1}t^{n-j+1}\over 1-q^{i}t^{n-j}}.
\end{equation}

\subsection{Skew symmetric functions}

Let us define as in \cite{Macdo} VI 7, the skew $Q$ functions by

\begin{equation}
\langle Q_{\lambda/\mu}, P_\nu\rangle_{q,t}
:={\langle Q_{\lambda},
P_\mu P_\nu\rangle_{q,t}}.
\end{equation}
Straightforwardly, one has
\begin{equation}\label{SQ2Q}
Q_{\lambda/\mu}(\X;q,t)=\sum_\nu\langle
Q_\lambda,P_\nu P_{\mu}\rangle_{q,t}
Q_\nu(\X;q,t).
\end{equation}
And classically, the following property hold. \footnote{See {\it e.g.} \cite{Macdo} VI.7 for a short
proof of this identity}
\begin{proposition}
Let $\X$ and $\Y$ be two alphabets, one has
 \[Q_{\lambda}(\X+\Y;q,t)=\sum_\mu Q_\mu(\X;q,t)Q_{\lambda/\mu}(\Y;q,t),\]
 or equivalently
\[P_{\lambda}(\X+\Y;q,t)=\sum_\mu P_\mu(\X;q,t)P_{\lambda/\mu}(\Y;q,t).\]
\end{proposition}
Equalities (\ref{Kqt}) and (\ref{lambda_t}) are generalized by identities \ref{skeqKqt} and
\ref{skewlambdaQ} as shown in \cite{Macdo}
 example 6 p.352
 \begin{equation}\label{skeqKqt}
\sum_\rho P_{\rho/\lambda}(\X;q,t)Q_{\rho/\mu}(\Y;q,t)=K_{qt}(\X,\Y)\sum_\rho P_{\mu/\rho}(\X;q,t)
Q_{\lambda/\rho}(\Y;q,t),
 \end{equation}
\begin{equation}\label{skewlambdaQ}
\sum_{\rho}Q_{\rho'/\lambda'}(\X;t,q)Q_{\rho/\mu}(\Y;q,t)=\lambda_1(\X\Y)\sum_\rho Q_{\mu'/\rho'}(\X,t,q)
Q_{\lambda/\rho}(\Y;q,t).
\end{equation}

\section{The substitution $x^p\rightarrow S^p(\Y)$ and the Macdonald polynomials\label{Ssubsti}}

Let $\X=\{x_1,\dots,x_n\}$ be a finite alphabet and $\Y$ be an other
(potentially infinite) alphabet. For simplicity we will denote by
$\int_\Y$ the substitution
 \begin{equation}
 \int_\Y x^p = S^p(\Y),
 \end{equation}
for each $x\in \X$ and each $p\in\Z$.
\subsection{Substitution formula}
Let us define the symmetric function
\begin{equation}\label{defH}
{\cal H}^{n,k}_{\lambda/\mu}(\Y;q,t):=\frac1{n!}\int_\Y
P_\lambda(\X;q,t)Q_\mu(\X^\vee;q,t)\Delta(\X,q,t)
\end{equation}
where $\X^\vee=\{x_1^{-1},\dots,x_n^{-1}\}$.

Set $\Y^{tq}:={1-t\over 1-q}\Y$ and consider the substitution
\begin{equation}
\int_{\Y^{tq}}x^p=S^p\left(\Y^{tq}\right)=Q_p(\Y;q,t).
\end{equation}
One has the following property.
%
%
\begin{theorem}\label{TqtSJT}
Let $\X=\{x_1,\dots,x_n\}$ be an alphabet, $\lambda=(\lambda_1,\dots,\lambda_n)$ be
a partition and $\mu\subset\lambda$. The polynomial ${\cal
H}_{\lambda/\mu}^{n,k}(\Y^{tq};q,t)$ is the Macdonald polynomial
\begin{equation}\label{qtSJT}
{\cal
H}_{\lambda/\mu}^{n,k}(\Y^{tq};q,t)=\frac1{n!}\prod_{(i,j)\in\lambda}{1-q^{i-1}t^{n-j+1}\over
1-q^it^{n-j}}\CT\{\Delta(\X,q,t)\}Q_{\lambda/\mu}(\Y,q,t)
\end{equation}
\end{theorem}
{\bf Proof}
From the definition of the $Q_\lambda$, one has
\begin{equation}
\int_{\Y^{tq}}x^p=Q_p(\Y;q,t)=\CT\{x^{-p}K_{q,t}(x,\Y)\}.
\end{equation}
 Hence, the polynomial ${\cal
H}_{\lambda/\mu}^{n,k}(\Y^{qt},q,t)$ is the constant term
\[
{\cal
H}_{\lambda/\mu}^{n,k}(\Y^{tq};q,t)=\frac1{n!}\CT\{P_{\lambda}(\X^\vee;q,t)
Q_\mu(\X;q,t)K_{q,t}(\X,\Y)\Delta(\X,q,t)\}.
\]
 As a special case of Equality (\ref{skeqKqt}),
\[
K_{qt}(\X,\Y)Q_{\mu}(\X;q,t)=\sum_\rho
P_{\rho/\mu}(\Y;q,t)Q_\rho(\X;q,t),
\]
holds and implies
\begin{equation}\label{I1}
\begin{array}{rcl}
{\cal
H}_{\lambda/\mu}^{n,k}(\Y^{tq},q,t)&=&\displaystyle\langle P_{\lambda}(\X^\vee;q,t), \sum_\rho
P_{\rho/\mu}(\Y;q,t)Q_\rho(\X;q,t) \rangle'_{q,t;n}\\
&=& \displaystyle\langle P_{\lambda}(\X^\vee;q,t), Q_\lambda(\X;q,t)\rangle_{q,t;n}'Q_{\lambda/\mu}(\Y,q,t).
\end{array}
\end{equation}
Equality (\ref{PQprime}) ends the proof. $\Box$\\ \\
  %
  %

%
%
%
%
 %

\subsection{Substitution dual formula}

Setting  $\overline\Y=\{-y_1,\dots,-y_m,\dots\}$ if $\Y=\{y_1,\dots,y_m,\dots\}$
\footnote{The operation $\Y\rightarrow \overline\Y$ makes sense for virtual alphabet since it sends
any homogeneous symmetric polynomial $P(\Y)$ of degree $p$ to $(-1)^pP(\Y)$.}, one observes the following propery.

%
%
\begin{theorem}
Let $\X=\{x_1,\dots,x_n\}$ be an alphabet, $\lambda=(\lambda_1,\dots,\lambda_n)$ be a partition and $\mu\subset\lambda$. One has
\begin{equation}\label{qtSdJT}
{\cal
H}_{\lambda/\mu}^{n,k}(-\overline\Y;q,t)={\cal
H}_{\lambda'/\mu'}^{n,k}(\Y^{qt};t,q)
\end{equation}
where $\Y^{qt}={1-q\over 1-t}\Y$.
\end{theorem}
{\bf Proof}
It suffices to show that
\[{\cal
H}_{\lambda/\mu}^{n,k}(-\overline\Y;q,t)=\frac1{n!}\prod_{(i,j)\in\lambda}{1-q^{i-1}t^{n-j+1}\over
1-q^it^{n-j}}\CT\{\Delta(\X,q,t)\}Q_{\lambda'/\mu'}(\Y,t,q).\]
The proof of this identity is almost the same than the proof of (\ref{qtSJT}) except than one uses the formula
\[
\prod(1+x_iy_j)Q_{\mu}(\X;t,q)=\sum_\rho
Q_\rho(\X;q,t)Q_{\rho'/\mu'}(\Y;t,q),
\]
which is a special case of identity (\ref{skewlambdaQ}).
$\Box$\\ \\
  %
  %

Note that in the case of partitions, one has
\begin{corollary}\label{CqtdJT}
\begin{equation}\label{qtdJT}
{\cal
H}_{\lambda}^{n,k}(-\overline\Y,q,t)=\frac1{n!}\prod_{(i,j)\in\lambda}{1-q^{i-1}t^{n-j+1}\over
1-q^it^{n-j}}\CT\{\Delta(\X,q,t)\}Q_{\lambda'}(\Y,t,q)
\end{equation}
\end{corollary}
 %
\begin{example}
\rm Consider the following equality
$$
{\cal H}^{2,3}_{41/3}(-\overline\Y;q,t)=(*)\CT\{\Delta(\X,q,t)\}Q_{2111/111}(\Y;t,q).
$$
where $\X=\{x_1,x_2\}$.
The coefficient $(*)$ is computed as follows. One writes the partition $[41]$ in a rectangle
of height $2$ and length $4$.
\[
\begin{array}{|c|c|c|c|}
\hline\times&\ &\ &\ \\
\hline\times&\times&\times&\times\\\hline
\end{array}
\]
Each $\times$ of coordinate $(i,j)$ is read as
the fraction $[i,j]:={1-q^{i-1}t^{3-j}\over
1-q^it^{2-j}}$. Hence
$$(*)=[1,2][1,1][2,1][3,1][4,1]={(1-t)(1-t^2)(1-qt^2)(1-q^2t^2)(1-q^3t^2)\over
(1-q)(1-qt)(1-q^2t)(1-q^3t)(1-q^4t)
}$$
\end{example}
\section{A formula relying the polynomials $P_{\lambda}\left({1-q\over 1-q^k}\X;q,q^k\right)$ and
$P_{\lambda}\left(\X;q,q^k\right)$\label{Smain}}
%
%
When $t=q^k$ with $k\in\N$, Corollary \ref{CqtdJT} gives
\begin{corollary}\label{CqqkJT}
\begin{equation}\label{qqkJT}
{\cal
H}_{\lambda}^{n,k}(-\overline\Y,q,q^k)=\beta_\lambda^{n,k}(q)Q_{\lambda'}(\Y;q^k,q).
\end{equation}
where
\[\beta_\lambda^{n,k}(q)=\prod_{i=0}^{n-1}\left[\lambda_{n-i}-1+k(i+1)\atop k-1\right]_q\]
 and $\left[n\atop p\right]_q={(1-q^n)\dots (1-q^{n-p+1})\over (1-q)\dots(1-q^r)}$ denotes the $q$-binomial.
\end{corollary}
{\bf Proof} From Corollary \ref{CqtdJT}, it remains to compute
$\CT\{\Delta(\X,q,t)\}$. The evaluation of this term is deduced from the $q$-Dyson conjecture \footnote{see \cite{Z1}
for a proof.}
\[
\CT\{\Delta(x;q,q^k)\}=n!\prod_{i=1}^n\left[ik-1\atop k-1\right]_q,
\]
and can be found  in \cite{Macdo} examples 1 p 372.

Hence,
\[{\cal
H}_{\lambda}^{n,k}(-\overline\Y,q,q^k)=\beta_\lambda^{n,k}(q)Q_{\lambda'}(\Y,q^k,q),\]
where
\begin{equation}\label{betaq}
 \beta_\lambda^{n,k}(q)=
\prod_{(i,j)\in\lambda}{1-q^{i+k(n-j+1)-1}\over
1-q^{i+k(n-j)}}\prod_{i=1}^n\left[ik-1\atop k-1\right]_q\\.
\end{equation}
But,
\[
\prod_{(i,j)\in\lambda}{1-q^{i+k(n-j+1)-1}\over 1-q^{i+k(n-j)}}=
\prod_{i=0}^{n-1}\prod_{j=1}^{\lambda_{n-i}}{1-q^{j+k(i+1)-1}\over
1-q^{j+ki}}.
\]
Hence, rearranging the factors appearing in the right hand side of
Equality (\ref{betaq}), one obtains
\begin{equation}\begin{array}{rcl}
 \beta_\lambda^{n,k}(q)&=&\displaystyle
\prod_{i=0}^{n-1}\left(\displaystyle\left[(i+1)k-1\atop
ik\right]_q\prod_{j=1}^{\lambda_{n-i}}{1-q^{j+k(i+1)-1}\over
1-q^{j+ki}}\right)\\
&=& \displaystyle \prod_{i=0}^{n-1}\left[\lambda_{n-i}-1+k(i+1)\atop
k-1\right]_q.\end{array}
 \end{equation}
 This ends the proof.$\Box$ 

\begin{example}
\rm
Set $k=2, n=3$ and consider the polynomial
\[
{\cal H}^{3,2}_{[320]}(-\overline\Y;q,q^2)=
\frac1{n!}\int_{-\overline\Y}P_{[32]}(x_1+x_2+x_3;q,q^2)\prod_{i\neq j}(1-x_ix_j^{-1})(1-qx_ix_j^{-1}).
\]
One has,
\[
{\cal H}^{3,2}_{[320]}(-\overline\Y;q,q^2)=
{\frac {\left (1-{q}^{5}\right )\left (1-{q}^{8}\right )}{\left (1-q
\right )^{2}}}
Q_{[221]}(\Y;q^2,q).
\]
\end{example}
 Let
 \begin{equation}\label{defOmega}
 \Omega_S:=\frac1{n!}\int_\X\prod_{i\neq j}(1-{x_ix_j^{-1}})
 \end{equation}
  and for each
 $v\in\Z^n$,\[ \tilde S_v(\X)=\det\left(x_i^{v_j+n-j}\right)\prod_{i<j}(x_i-x_j)^{-1}.\]
\begin{lemma}\label{LOmega}
 If $v$ is  any vector in $\Z^n$, one has
 \begin{equation}\label{Omega}
\Omega_S \tilde S_v(\X)=
S_v(\X):=\det(S^{v_i-i+j}(\X))
 \end{equation}
\end{lemma}
{\bf Proof}
The identity is obtain by the direct computation:
\[\begin{array}{rcl}
\displaystyle\frac1{n!}\int_\X \tilde S_v(\X)\prod_{i<j}\left(1-x_ix_j^{-1}\right)&=&\displaystyle\frac1{n!}
\int_\X\det\left(x_i^{v_j-j+1}\right)\det(x_i^{j-1})\\
&=&\displaystyle\frac1{n!}\int_\X\sum_{\sigma_1,\sigma_2\in\S_n}{\rm sign}(\sigma_1\sigma_2)\prod_{i}x^{v_{\sigma_1(i)}
-\sigma_1(i)+\sigma_2(i)-1}\\
&=&\displaystyle \frac1{n!}\sum_{\sigma_1\sigma_2}{\rm sign}(\sigma_1\sigma_2)\prod_iS_{v_{\sigma_1(i)}-\sigma_1(i)+\sigma_2(j)}\\
&=&\det(S^{v_i-i+j}(\X)).
\end{array}
\]
$\Box$\\ \\
In particular, $\Omega_S$ lets invariant any symmetric polynomial.
The operator
\begin{equation}
\label{defA}{\goth A}_m:=\Omega_S\Lambda^n(\X)^{-m}
\end{equation}
acts on symmetric  polynomials by  substracting
$m$ on each part of partitions appearing in their expansion in the Schur basis.
\begin{example}\rm
If $\X=\{x_1,x_2,x_3\}$ consider the polynomial, and $\lambda=[320]$. One has

$$P_{32}(\X;q,t)=S_{{32}}(\X)+{\frac {\left (-q+t\right )S_{{311}}(\X)}{qt-1}}+{\frac {
\left (q+1\right )\left (q{t}^{2}-1\right )\left (-q+t\right )S_{{22
1}}(\X)}{\left (qt-1\right )^{2}\left (qt+1\right )}}.
$$
Hence,
$$\begin{array}{rcl}
{\goth A}_1P_{32}(\X;q,t) &=&
                              {\frac {\left (-q+t\right )S_{{2}}(\X)}{qt-1}}+{\frac {\left (q+1
                               \right )\left (q{t}^{2}-1\right )\left (-q+t\right )S_{{11}}(\X)}{
                               \left (qt-1\right )^{2}\left (qt+1\right )}}\\
                         &=&  {\frac {\left (-q+t\right )\left (t+1\right )\left ({q}^{2}t-1\right )
                               P_{{11}}(\X;q,t)}{\left (qt-1\right )^{2}\left (qt+1\right )}}+{\frac {\left
                               (-q+t\right )P_{{2}}(\X;q,t)}{qt-1}}.
\end{array}
$$

\end{example}

\begin{theorem}\label{th-X}
If $\lambda$ denotes a partition of length at most $n$, one has
\begin{equation}
{\cal A}_{(k-1)(n-1)}P_\lambda(\X;q,q^k)\prod_{l=1}^{k-1}\prod_{i\neq j}(x_i-q^lx_j)=
\beta_{\lambda}^{n,k}(q)P_{\lambda}\left({1-q\over 1-q^k}\X;q,q^k\right)
\end{equation}
\end{theorem}
{\bf Proof}
From the definitions of the operators ${\cal A}_m$ (\ref{defA}) and $\Omega_S$ (\ref{defOmega}), one obtains
\[
{\cal A}_{(k-1)(n-1)}P_\lambda(\X;q,q^k)\prod_{l=1}^{k-1}\prod_{i\neq j}(x_i-q^lx_j)
=\frac1{n!}\int_{\X} P_{\lambda}(\X;q,q^k)\Delta(\X,q,q^k).
\]
 Corollary \ref{CqqkJT} implies
\[\frac1{n!}\begin{array}{rcl}
\displaystyle\int_{\X} P_{\lambda}(\X;q,q^k)\Delta(\X,q,q^k)&=&{\cal H}_{\lambda}^{n,k}(\X;q,q^k)
\\&=& \displaystyle\beta_\lambda^{n,k}(q)Q_{\lambda'}(-\overline\X;q^k,q)\end{array}
\]
But, from Lemma \ref{P-X2Pqt}, one has
\[Q_{\lambda'}(-\overline\X;q^k,q)=(-1)^{|\lambda|}Q_{\lambda'}(-\X;q^k,q)=P_{\lambda}\left({1-q\over 1-q^k}\X;q,q^k\right).\]

The result follows.
$\Box$
  %
  %
\begin{example}\rm
Set $k=2$, $n=3$ and $\lambda=[2]$. One has
\[\begin{array}{l}
\displaystyle P_{[2]}(x_1+x_2+x_3;q,q^2)\prod_{i\neq j}(x_i-qx_j)=\displaystyle -q^3S_{[6,2]}
+q^2\frac{q^3-1}{q-1}S_{[6,1,1]}
\\\displaystyle+\frac{q^2(q^5-1)}{q^3-1}S_{[5,3]} -
\frac{q(q^2+1)(q^5-1)}{q^3-1}S_{[5,2,1]}-\frac{q(q^7-1)}{q^3-1}S_{[4,3,1]}
 +\frac{q^7-1}{q-1}S_{[4,2,2]}.\end{array}
\]
And,
\[
{\cal A}_{2}P_{[2]}(x_1+x_2+x_3;q,q^2)\prod_{i\neq j}(x_i-qx_j)=\frac{q^7-1}{q-1}S_{[2]}.
\]
Since,
\[
P_{[2]}\left(\frac{x_1+x_2+x_3}{1+q};q,q^2\right)=\frac {1-q}{1-q^3}S_{[2]}
\]
one obtains
\[
{\cal A}_{2}P_{[2]}(x_1+x_2+x_3;q,q^2)\prod_{i\neq j}(x_i-qx_j)=\left[1\atop 1\right]_q
\left[3\atop 1\right]_q\left[7\atop 1\right]_qP_{[2]}\left(\frac{x_1+x_2+x_3}{1+q};q,q^2\right).
\]

\end{example}

As   a consequence, one has
\begin{corollary}\label{mu+rect}
If $\lambda=\mu+[((k-1)(n-1))^n]$,
\[
P_\mu(\X;q,q^k)\prod_{l=1}^{k-1}\prod_{i\neq j}(x_i-q^lx_j)=
\beta_{\lambda}^{n,k}(q)P_{\lambda}\left({1-q\over 1-q^k}\X;q,q^k\right)
.\]
\end{corollary}
{\bf Proof} Since the size of $\X$ is $n$,
\[P_{\lambda}(\X;q,q^k)=P_\mu(\X;q,q^k)(x_1\dots x_n)^{(k-1)(n-1)}.\]
Then, the result is a direct consequence of Theorem \ref{th-X}.$\Box$

\begin{example}\rm
Set $k=3$, $n=2$ and $\lambda=[5,2]$. One has,
\[
\begin{array}{l}
\displaystyle P_{[5,2]}(x_1+x_2;q,q^3)(x_1-qx_2)(x_1-q^2x_2)(x_2-qx_1)(x_2-q^2x_1)=\\
\displaystyle {q}^{3}S_{[9,2]}+\frac{(1-q^{7})(1+q^4)}{1-q^5} S_{[7,4]}-
\frac{(1-q^2)(1+q)(1+q^2)(1+q^4)}{1-q^5 }S_{[8,3]}.
\end{array}
\]
This implies
\[\begin{array}{l}\displaystyle {\cal A}_2P_{[5,2]}(x_1+x_2;q,q^3)(x_1-qx_2)(x_1-q^2x_2)(x_2-qx_1)(x_2-q^2x_1)=
\\\displaystyle (x_1x_2)^{-2}P_{[5,2]}(x_1+x_2;q,q^3)(x_1-qx_2)(x_1-q^2x_2)(x_2-qx_1)(x_2-q^2x_1)=
\\\displaystyle P_{[3]}(x_1+x_2;q,q^3)(x_1-qx_2)(x_1-q^2x_2)(x_2-qx_1)(x_2-q^2x_1).\end{array}\]
One verifies that
\[\begin{array}{l}\displaystyle
P_{[3]}(x_1+x_2;q,q^3)(x_1-qx_2)(x_1-q^2x_2)(x_2-qx_1)(x_2-q^2x_1)=\\\displaystyle
\left[4\atop 2\right]_q\left[10\atop 2\right]_qP_{[5,2]}(\frac{x_1+x_2}{1+q+q^2};q,q^3).\end{array}
\]
\end{example}

\begin{remark}\rm
If $\mu$ is the empty partition, Corollary \ref{mu+rect} gives
\begin{equation}\nonumber\label{qdisc}
\prod_{l=1}^{k-1}\prod_{i\neq j}(x_i-q^lx_j)=
\beta_{\lambda}^{n,k}(q)P_{[((k-1)(n-1))^n]}\left({1-q\over 1-q^k}\X;q,q^k\right).
\end{equation}
This equality generalizes an identity given in \cite{BBL}:
$$\prod_{i< j}(x_i-x_j)^{2(k-1)}={(-1)^{((k-1)n(n-1)\over2}\over n!}
\left(kn\atop k,\dots,k\right)P_{n^{(n-1)(k-1)}}^{(k)}(-\X),$$
where $P_\lambda^{(k)}(\X)=\displaystyle\lim_{q\rightarrow1} P^{(\alpha)}_\lambda(\X;q,q^k)$ denotes a Jack polynomial
(see {\it e.g.} \cite{Macdo}).\\
The expansion of the powers of the discriminant and their $q$-deformations in  different basis of symmetric
functions is a difficult problem having many applications, for example, in the study of Hua-type integrals (see {\it e.g.}
\cite{Kan,Kor}) or
in the context of the factional quantum Hall effect ({\it e.g.} \cite{DGIL,KTW,Lau,STW}).\\
Note that in \cite{BL}, we gave an expression of an other $q$-deformation of the powers of the discriminant
as  staircase
Macdonald polynomials. This deformation is also relevant in the study of the expansion of $\prod_{i< j}(x_i-x_j)^{2k}$ in the
Schur basis, since we generalized  \cite{BL} a result of \cite{KTW}.
\end{remark}

\section{Macdonald polynomials at $t=q^k$ as eigenfunctions\label{Seigen}}

Let $\Y=\{y_1,\dots,y_{kn}\}$ be an alphabet of cardinality $kn$ with $y_1=x_1,\ \dots,\ y_n=x_n$.
 One considers the symmetrizer  $\pi_\omega $
defined by
\[
\pi_\omega f(y_1,\dots,y_{kn})=
\prod_{i<j}(x_i-x_j)^{-1}\sum_{\sigma\in\S_{kn}}\sign(\sigma)f(y_{\sigma(1)},\dots,y_{\sigma(kn)})
y_{\sigma(1)}^{kn-1}\dots y_{\sigma(kn-1)}.
\]
Note that $\pi_\omega$ is the isobaric divided difference associated to the maximal permutation $\omega$ in
$\S_{kn}$.


This operator applied to a  symmetric function of the alphabet $\X$
increases the alphabet from $\X$ to $\Y$ in its expansion in the
Schur basis, since
\begin{equation} \label{X2Y}
\pi_\omega S_\lambda(\X)=S_\lambda(\Y).
\end{equation}
Indeed,  the image of the monomial $y_{1}^{i_1}\dots y_{kn}^{i_{kn}}$
is the Schur function $S_I(\Y)$. Since
\[
\pi_\omega S_\lambda(\X)=\pi_\omega x_1^{\lambda_1}\dots x_{n}^{\lambda_{n}}=
\pi_\omega y_1^{\lambda_1}\dots y_n^{\lambda_n} y_{n+1}^0\dots y_{kn}^0,
\]
one recovers Equality (\ref{X2Y}).

One defines the operator $\pi^{tq}$ which consists in applying $\pi_\omega$ and specializing the result to the alphabet
$$\X^{tq}:=\{x_1,\dots,x_n, qx_1,\dots,qx_n,\dots,q^{k-1}x_1,\dots,q^{k-1}x_n\}.$$
From Equality (\ref{X2Y}), one has
\begin{equation}\label{S2Sqqk}
\pi^{tq}_\omega S_\lambda(\X)=S_\lambda\left((1+q+\dots+q^{k-1})\X\right),
\end{equation}
for $l(\lambda)\leq n$.
Furthermore, the expansion of $S_\lambda\left((1+q+\dots+q^{k-1})\X\right)$
in the Schur basis being triangular, the operator $\pi^{tq}$ defines an automorphism of the
space $Sym_{\leq n}$ generated by the Schur functions indexed by partitions whose length
are less or equal to $n$, {\it i.e.} for each function $f\in Sym_{\leq n}$, one has
\begin{equation}\label{X2X^tq}
\pi^{tq} f(\X)=f(\X^{tq}).
\end{equation}
In particular, one has

\begin{lemma}\label{Pnabla}
 Let $\lambda$ be a partition such that $l(\lambda)\leq
n$ then
\begin{equation} \pi_\omega^{tq} P_\lambda\left({1-q\over
1-q^k}\X;q,t=q^k\right)=P_\lambda(\X,q,q^k).
\end{equation}
\end{lemma}
{\bf Proof}
It suffices to remark that
$P_\lambda\left({1-q\over
1-q^k}\X;q,q^k\right) \in Sym_{\leq n}(\X)$.\footnote{This can be seen as a consequence of the
determinantal expression of the expansion of
 $P_\lambda(\X,q,t)$ in the Schur basis evaluated on the alphabet $\X^{tq}$ (see \cite{LLM}).}\\
 It follows from (\ref{X2X^tq}),
\[\pi^{tq}_\omega P_\lambda\left({1-q\over
1-q^k}\X;q,q^k\right)=P_\lambda\left({1-q\over
1-q^k}\X^{tq};q,q^k\right)=P_\lambda(\X,q,q^k).\]
 $\Box$\\ \\
Consider the operator ${\cal M}:f\rightarrow {\cal M}f$ defined by
\[
{\cal M}:=(x_1\dots x_n)^{(k-1)(1-n)}\pi_\omega^{tq}\prod_{l=1}^{k-1}\prod_{i\neq j}(x_i-q^lx_j).
\]
The following theorem shows that the Macdonald polynomials are the eigenfunctions of the operator ${\cal M}$.
\begin{theorem}\label{eigen1}
The Macdonald polynomials $P_\lambda(\X;q,q^k)$ are eigenfunctions of $\cal M$.
The eigenvalue associated to $P_\mu(\X;q,q^k)$ is $\beta_{\mu+((k-1)(n-1))^n}^{n,k}(q)$.
Furthermore, if $k>1$, the dimension of
each eigenspace is $1$.
\end{theorem}
{\bf Proof}
From Corollary \ref{mu+rect}, one has
\[
P_\mu(\X;q,q^k)\prod_{l=    1}^{k-1}\prod_{i\neq j}(x_i-q^lx_j)=
\beta_{\lambda}^{n,k}(q)P_{\lambda}\left({1-q\over 1-q^k}\X;q,q^k\right)
\]
where $\lambda=\mu+((k-1)(n-1))^n$.
Applying $\pi_\omega^{tq}$ to the left and the right hand sides of this equality, one obtains
from Lemma \ref{Pnabla}
\[\begin{array}{rcl} \displaystyle
\pi_\omega^{tq} P_\mu(\X;q,q^k)\prod_{l=1}^{k-1}\prod_{i\neq j}(x_i-q^lx_j)&=&
\beta_{\lambda}^{n,k}(q)\pi_\omega^{tq} P_{\lambda}\left({1-q\over 1-q^k}\X;q,q^k\right)\\
&
=&\beta_{\lambda}^{n,k}(q) P_{\lambda}\left(\X;q,q^k\right).\end{array}
\]
Since the cardinality of $\X$ is $n$, one has
\[
P_{\lambda}\left(\X;q,q^k\right)=(x_1\dots x_n)^{(k-1)(n-1)}P_{\lambda}\left(\X;q,q^k\right),
\]
and
\begin{equation}\label{MP}
{\cal M}
P_\mu(\X;q,q^k)=\beta_{\mu+[((k-1)(n-1))^n]}^{n,k}(q)P_\mu(\X;q,q^k).
\end{equation}

Suppose now that $k> 1$. It remains to prove that the dimensions of
the eigenspaces equal $1$. More precisely, It suffices to show that
 $\beta_\lambda(q)=\beta_\mu(q)$ implies $\lambda=\mu$. The
 denominators of $\beta_\lambda(q)$ and $\beta_\mu(q)$ being the
 same, one needs  only to examine the numerators, that is the products
 $\gamma_\lambda=\prod_{i=0}^{n-1}(q^{\lambda_{n-i}+ki};q)_{k-1}$
 and $\gamma_\mu=\prod_{i=0}^{n-1}(q^{\mu_{n-i}+ki};q)_{k-1}$.
One needs the following lemma.
\begin{lemma}\label{1-q}
Let $I=\{i_1,\dots,i_n\}$ and $J=\{j_1,\dots,j_m\}$ be two finite subsets of $\N\setminus\{0\}$.
Then, $I\neq J$ implies $\prod_{i\in I}(1-q^i)\neq \prod_{j\in J}(1-q^j).$
\end{lemma}
{\bf Proof}
Without lost of generalities, one can suppose $I\cap J=\emptyset$. Suppose that $i_1\leq\dots\leq i_n$ and
 $j_1\leq\dots\leq j_m$. Then, expanding the two products, one finds
 \[
\prod_{i\in I}(1-q^i)=1-q^{i_1}+\sum_{l>i_1} (*)q^{l}\neq 1-q^{j_1}+\sum_{l>j_1} (*) q^{l}=\prod_{j\in J}(1-q^j).
 \]
$\Box$\\
 Each
 term $(q^{\lambda_{n-i}+ki};q)_{k-1}$ is characterized by the
 degree of its factor of lower degree : $\lambda_{n-i}+ki$. Hence, from Lemma \ref{1-q},
 $\beta_\lambda(q)=\beta_\mu(q)$ implies that it exists a
 permutation $\sigma$ of $\S_n$ verifying
 \[
\lambda_i+k(n-i)=\mu_{\sigma_i}+k(n-\sigma_i),
 \]
for each $i$. But, since $\lambda$ is decreasing, one has
\[
\lambda_i+k(n-i)-\lambda_{i-1}-k(n-i+1)\leq 0.
\]
And then,
\begin{equation}\label{mumu}
\mu_{\sigma_i}+k(n-\sigma_i)-\mu_{\sigma_{i-1}}-k(n-\sigma_{i-1})\leq
0.
\end{equation}
But, since $\mu$ is decreasing, $\sigma_{i-1}-\sigma_i$ has the same
sign than $\mu_{\sigma_i}-\mu_{\sigma_{i-1}}$. As a consequence,
Inequality (\ref{mumu}) implies $\sigma_i>\sigma_{i-1}$ for each
$i$. The only possibility is $\sigma=Id$, which ends the
proof.$\Box$
\begin{example}
\rm
If $n=5$, the eigenvalues associated to the partitions of $4$ are
\[
\begin{array}{l}
\beta_{[4\,k,4\,k-4,4\,k-4,4\,k-4,4\,k-4]}^{4,k}=\left[5k-5\atop k-1\right]_q
\left[6k-5\atop k-1\right]_q\left[7k-5\atop k-1\right]_q\left[8k-5\atop k-1\right]_q\left[9k-1\atop k-1\right]_q,\\
\beta_{[4\,k-1,4\,k-3,4\,k-4,4\,k-4,4\,k-4]}^{4,k}=\left[5k-5\atop k-1\right]_q\left[6k-5\atop k-1\right]_q
\left[7k-5\atop k-1\right]_q\left[8k-4\atop k-1\right]_q\left[9k-2\atop k-1\right]_q,\\
\beta_{[4\,k-2,4\,k-2,4\,k-4,4\,k-4,4\,k-4]}^{4,k}=\left[5k-5\atop k-1\right]_q
\left[6k-5\atop k-1\right]_q\left[7k-5\atop k-1\right]_q\left[8k-3\atop k-1\right]_q\left[9k-3\atop k-1\right]_q,\\
\beta_{[4\,k-2,4\,k-3,4\,k-3,4\,k-4,4\,k-4]}^{4,k}=\left[5k-5\atop k-1\right]_q\left[6k-5\atop k-1\right]_q
\left[7k-4\atop k-1\right]_q\left[8k-4\atop k-1\right]_q\left[9k-3\atop k-1\right]_q,\\
\beta_{[4\,k-3,4\,k-3,4\,k-3,4\,k-3,4\,k-4]}^{4,k}=\left[5k-5\atop k-1\right]_q\left[6k-4\atop k-1\right]_q
\left[7k-4\atop k-1\right]_q\left[8k-4\atop k-1\right]_q\left[9k-4\atop k-1\right]_q.
\end{array}
\]
\end{example}
\section{Expression of $\cal M$ in terms of Cherednik elements\label{SChered}}
In this paragraph, we restate Proposition \ref{eigen1} in terms of Cherednik operators.
Cherednik's operators $\{\xi_i;i\in \{1,\dots,n\}\}=:\Xi$ are commutative
elements of the double affine Hecke algebra. The Macdonald
polynomials $P_\lambda(\X;q,t)$ are eigenfunctions of symmetric polynomials   $f(\Xi)$ and
the eigenvalues are obtained substituting each occurrence of $\xi_i$ in $f(\Xi)$ by $q^{\lambda_i}t^{n-i}$
 (see \cite{Lasc2} for more details).\\
Suppose that $k>1$ and
consider the operator $\tilde{\cal M}:f\rightarrow \tilde{\cal M}f$ defined by
\begin{equation}\label{tildeM2M}
\tilde{\cal M}:=\prod_{i=1}^{k-1}(1-q^i)^n{\cal M}.
\end{equation}
From Proposition \ref{eigen1}, one has
\begin{equation}
\tilde{\cal M}P_{\lambda}(\X;q,q^k)=\prod_{i=0}^{n-1}\prod_{j=1}^{k-1}(1-q^{\lambda_{n-i}+k(i+1)-j})P_{\lambda}(\X;q,q^k).
\end{equation}
The following proposition shows that $\tilde{\cal M}$ admits a closed expression in terms of Cherednick elements.
%
%
\begin{proposition}\label{Teigen} One supposes that $k>1$.
For any symmetric function $f$, one has
\begin{equation}\label{tM2Ch}
\tilde{\cal M}f(\X)=\prod_{l=1}^{k-1}\prod_{i=1}^n(1-q^{l+k}\xi_i)f(\X).
\end{equation}
\end{proposition}
{\bf Proof}
From Theorem \ref{eigen1}, it suffices to prove  the formula  (\ref{tM2Ch}) for $f=P_\lambda$.
 The polynomial $P_\lambda(\X;q,t)$ is an eigenfunction of the operator
$\prod_{l=1}^{k-1}\prod_{i=1}^n(1-q^{l+k}\xi_i)$ and its  eigenvalues is
$\prod_{l=1}^{k-1}\prod_{i=1}^n (1-q^{l+\lambda_i}t^{n-i+1})P_\lambda(\X;q,t)$.
Hence, setting $t=q^k$, we obtain
\[
\prod_{l=1}^{k-1}\prod_{i=1}^n(1-q^{l+k}\xi_i)P_\lambda(\X;q,q^k)=\prod_{l=1}^{k-1}\prod_{i=1}^n
(1-q^{l+\lambda_i+k(n-i+1)})P_\lambda(\X;q,q^k).
 \]
Comparing this expression to Equality (\ref{tildeM2M}), one finds the result.
$\Box$
  %
  %
\\ \\
\noindent{\bf Acknowledgments}
The author is  grateful to Alain Lascoux and Jean-Yves Thibon for  fruitful discussions.

 \end{document}